\documentclass[a4paper,11pt]{article}
\usepackage{amssymb}
\usepackage{amsfonts}
\usepackage{amsmath}
\usepackage{geometry}
\usepackage{hyperref}

\input{tcilatex}
\begin{document}

\title{\textbf{Relative Error Control in Bivariate Interpolatory Cubature}}
\author{J.S.C. Prentice \\
%EndAName
{\small (ORCID: 0000-0003-2506-7327)}\\
Senior Research Officer\\
Mathsophical Ltd.\\
Johannesburg, South Africa\\
{\small jpmsro@mathsophical.com}}
\maketitle

\begin{abstract}
We describe an algorithm for controlling the relative error in the numerical
evaluation of a bivariate integral, without prior knowledge of the magnitude
of the integral. In the event that the magnitude of the integral is less
than unity, absolute error control is preferred. The underlying quadrature
rule is positive-weight interpolatory and composite. Some numerical examples
demonstrate the algorithm.
\end{abstract}

\begin{flushleft}
\textbf{{\Large {Keywords}}}

\smallskip

Bivariate integral; Composite quadrature; Interpolatory quadrature;
Cubature; Relative error; Absolute error; Error control.

\medskip \bigskip

\textbf{{\Large {MSC 2010}}}

\smallskip

65D30; 65D32; 65G20
\end{flushleft}

\section{Introduction}

Let $G\left( x,y\right) $ be such that%
\begin{equation*}
G:%
%TCIMACRO{\U{211d} }%
%BeginExpansion
\mathbb{R}
%EndExpansion
^{2}\rightarrow 
%TCIMACRO{\U{211d} }%
%BeginExpansion
\mathbb{R}
%EndExpansion
,
\end{equation*}%
$G$ is Riemann integrable in both of its variables, and $G$ and all its
derivatives relevant to this study exist on every region of integration
considered therein.

In this paper, we consider the evaluation of%
\begin{equation*}
I\left[ G\left( x,y\right) \right] \equiv
\int\limits_{a}^{b}\int\limits_{l\left( x\right) }^{u\left( x\right)
}G\left( x,y\right) dydx
\end{equation*}%
using cubature based on composite interpolatory quadrature, such that%
\begin{equation}
\left\vert \frac{I\left[ G\left( x,y\right) \right] -Q_{C}\left[ G\left(
x,y\right) \right] }{I\left[ G\left( x,y\right) \right] }\right\vert
\leqslant \varepsilon ,  \label{rel err < eps}
\end{equation}%
where $Q_{C}\left[ G\left( x,y\right) \right] $ is the composite cubature of 
$G\left( x,y\right) $, $\varepsilon $ is a user-imposed tolerance, and an
estimate of $I\left[ G\left( x,y\right) \right] $ is \textit{not} known a
priori. In other words, we seek to control the relative error in the
cubature, without prior estimation of the integral. The problem is easily
understood with reference to (\ref{rel err < eps}): we have%
\begin{equation*}
\left\vert I\left[ G\left( x,y\right) \right] -Q_{C}\left[ G\left(
x,y\right) \right] \right\vert \leqslant \varepsilon \left\vert I\left[
G\left( x,y\right) \right] \right\vert ,
\end{equation*}%
so that $\varepsilon \left\vert I\left[ G\left( x,y\right) \right]
\right\vert $ is an absolute tolerance. In principle, interpolatory methods
readily admit absolute error control but, since $I\left[ G\left( x,y\right) %
\right] $ is not known, we cannot impose $\varepsilon \left\vert I\left[
G\left( x,y\right) \right] \right\vert $ as a tolerance. Controlling the
relative error is appropriate when dealing with integrals of large
magnitude; for such integrals, absolute error control can be very
inefficient. Again, however, this presents a problem, since the magnitude of 
$I\left[ G\left( x,y\right) \right] $ is not known, so we do not even know
whether absolute or relative error control should be applied.

The algorithm we present here is, in a sense, a `first-principles' method,
since it is based entirely on classical concepts relating to interpolatory
quadrature. Nevertheless, it is an extension of our previous work regarding
univariate numerical integration \cite{jscp}. The list of references \cite%
{bf}$-$\cite{ss} is our bibliography, and is drawn from the established
literature.

A note regarding terminology: from this point onwards we will use the term 
\textit{quadrature} to refer to the numerical approximation of a univariate
integral, and the term \textit{cubature} to refer to the numerical
approximation of a multivariate integral.

\section{The Algorithm}

We transform $\left[ a,b\right] $ to $\left[ 0,1\right] $ by means of%
\begin{equation}
x=\left( b-a\right) w+a\equiv m_{1}w+a,  \label{transform x}
\end{equation}%
where $x\in \left[ a,b\right] ,w\in \left[ 0,1\right] $ and $m_{1}$ has been
implicitly defined.

If%
\begin{align*}
l_{1}& \equiv \min_{\left[ a,b\right] }\left\{ l\left( x\right) ,u\left(
x\right) \right\} \\
u_{1}& \equiv \max_{\left[ a,b\right] }\left\{ l\left( x\right) ,u\left(
x\right) \right\}
\end{align*}%
then the transformation between $\left[ l_{1},u_{1}\right] $ and $\left[ 0,1%
\right] $ is given by%
\begin{equation}
y=\left( u_{1}-l_{1}\right) z+l_{1}\equiv m_{2}z+l_{1},  \label{transform y}
\end{equation}%
where $y\in \left[ l_{1},u_{1}\right] ,z\in \left[ 0,1\right] $ and $m_{2}$
has been implicitly defined.

As a result of these affine transformations,%
\begin{equation*}
\int\limits_{a}^{b}\int\limits_{l\left( x\right) }^{u\left( x\right)
}G\left( x,y\right) dydx=\int\limits_{0}^{1}\int\limits_{\widetilde{l}\left(
w\right) }^{\widetilde{u}\left( w\right) }\widetilde{G}\left( w,z\right)
m_{1}m_{2}dzdw,
\end{equation*}%
where%
\begin{align*}
\widetilde{G}\left( w,z\right) & \equiv G\left( m_{1}w+a,m_{2}z+l_{1}\right)
\\
\widetilde{u}\left( w\right) & \equiv \frac{u\left( m_{1}w+a\right) -l_{1}}{%
m_{2}} \\
\widetilde{l}\left( w\right) & \equiv \frac{l\left( m_{1}w+a\right) -l_{1}}{%
m_{2}}.
\end{align*}

We determine%
\begin{equation}
M\equiv \max \left\{ 1,\max_{\widetilde{R}}\left\vert \widetilde{G}\left(
w,z\right) m_{1}m_{2}\right\vert \right\} ,  \label{M}
\end{equation}%
where $\widetilde{R}$ is the domain of integration defined by the transforms
(\ref{transform x}) and (\ref{transform y}). Note that $\widetilde{R}%
\subseteq \left[ 0,1\right] \times \left[ 0,1\right] .$ Hence, we define%
\begin{equation*}
g\left( w,z\right) \equiv \frac{\widetilde{G}\left( w,z\right) m_{1}m_{2}}{M}%
.
\end{equation*}

Now,%
\begin{align*}
& \left\vert I\left[ g\left( w,z\right) \right] -Q_{C}\left[ g\left(
w,z\right) \right] \right\vert \leqslant \varepsilon \\
& \Rightarrow \left\vert MI\left[ g\left( w,z\right) \right] -MQ_{C}\left[
g\left( w,z\right) \right] \right\vert \leqslant M\varepsilon \\
& \Rightarrow \left\vert I\left[ \widetilde{G}\left( w,z\right) m_{1}m_{2}%
\right] -Q_{C}\left[ \widetilde{G}\left( w,z\right) m_{1}m_{2}\right]
\right\vert \leqslant M\varepsilon \\
& \Rightarrow \left\vert I\left[ \widetilde{G}\left( w,z\right) m_{1}m_{2}%
\right] -Q_{C}\left[ \widetilde{G}\left( w,z\right) m_{1}m_{2}\right]
\right\vert \leqslant \left\vert \frac{I\left[ \widetilde{G}\left(
w,z\right) m_{1}m_{2}\right] }{I\left[ g\left( w,z\right) \right] }%
\right\vert \varepsilon \\
& \Rightarrow \frac{\left\vert I\left[ \widetilde{G}\left( w,z\right)
m_{1}m_{2}\right] -Q_{C}\left[ \widetilde{G}\left( w,z\right) m_{1}m_{2}%
\right] \right\vert }{I\left[ \widetilde{G}\left( w,z\right) m_{1}m_{2}%
\right] }\leqslant \frac{\varepsilon }{\left\vert I\left[ g\left( w,z\right) %
\right] \right\vert } \\
& \Rightarrow \frac{\left\vert I\left[ G\left( x,y\right) \right] -Q_{C}%
\left[ G\left( x,y\right) \right] \right\vert }{I\left[ G\left( x,y\right) %
\right] }\leqslant \frac{\varepsilon }{\left\vert I\left[ g\left( w,z\right) %
\right] \right\vert }\approx \frac{\varepsilon }{\left\vert Q_{C}\left[
g\left( w,z\right) \right] \right\vert }.
\end{align*}%
In the last inequality, we use the fact that changes in variable preserve
the value of both the integral and the quadrature-based cubature.

Clearly, from the last inequality, 
\begin{equation*}
\frac{\varepsilon }{\left\vert Q_{C}\left[ g\left( w,z\right) \right]
\right\vert }
\end{equation*}%
is an estimated bound on the relative error in $Q_{C}\left[ G\left(
x,y\right) \right] =MQ_{C}\left[ g\left( w,z\right) \right] .$ This estimate
is good if $Q_{C}\left[ g\left( w,z\right) \right] $ is accurate which, in
turn, is determined by the choice of $\varepsilon .$

Now, assuming $\left\vert I\left[ G\left( x,y\right) \right] \right\vert
\geqslant 1,$ 
\begin{equation*}
\frac{1}{\left\vert Q_{C}\left[ g\left( w,z\right) \right] \right\vert }=%
\frac{M}{\left\vert Q_{C}\left[ G\left( x,y\right) \right] \right\vert }%
\approx \frac{M}{\left\vert I\left[ G\left( x,y\right) \right] \right\vert }
\end{equation*}%
and, since $M$ is the maximum possible value of $\left\vert I\left[ G\left(
x,y\right) \right] \right\vert $ (by construction, see (\ref{M})), we have%
\begin{equation*}
\frac{\varepsilon }{\left\vert Q_{C}\left[ g\left( w,z\right) \right]
\right\vert }\sim \varepsilon ,
\end{equation*}%
provided $\left\vert I\left[ G\left( x,y\right) \right] \right\vert $ is not
substantially smaller than $M$. For many practical situations, this will be
the case. However, we have no prior knowledge of $\left\vert Q_{C}\left[
g\left( w,z\right) \right] \right\vert \approx \left\vert I\left[ g\left(
w,z\right) \right] \right\vert ,$ so we must be willing to accept the
estimate, whatever it may be. Obviously, we cannot expect that the relative
error will satisfy the tolerance $\varepsilon ,$ even if $\left\vert I\left[
g\left( w,z\right) \right] -Q_{C}\left[ g\left( w,z\right) \right]
\right\vert $ does. Note that if 
\begin{equation*}
\left\vert I\left[ g\left( w,z\right) \right] -Q_{C}\left[ g\left(
w,z\right) \right] \right\vert =\varepsilon ,
\end{equation*}
then $\frac{\varepsilon }{\left\vert Q_{C}\left[ g\left( w,z\right) \right]
\right\vert }$ is not merely an upper bound, but is a very good estimate of
the relative error itself.

If $\left\vert I\left[ G\left( x,y\right) \right] \right\vert <1,$ then the
relative error could be considerably larger than $\varepsilon ,$
particularly if $\left\vert I\left[ G\left( x,y\right) \right] \right\vert
\sim 0,$ but in this case we favour absolute error control (for reasons to
be discussed later), and so the relative error is not relevant. The quantity 
$M\varepsilon $ is an upper bound on the absolute error.

If the estimate of the absolute or relative error is considered too large,
say by a factor of $\eta ,$ then we simply redo the calculation, this time
with a tolerance of 
\begin{equation*}
\frac{\varepsilon }{\eta }.
\end{equation*}%
This refinement is a very important feature of the algorithm, since it
enables a desired tolerance to be achieved in a controlled manner, even if
it requires a repetition of the calculation. We are sure that such
repetition is a small price to pay for a solution of acceptable quality.

\section{Bivariate composite interpolatory cubature}

Here, we briefly describe bivariate composite interpolatory cubature,
including the relevant error analysis. We will consider the effect of
roundoff error on error control, and offer a criterion for choosing between
absolute and relative error control. A reasonable degree of familiarity with
interpolatory quadrature is assumed.

\subsection{The form of bivariate composite interpolatory cubature}

The composite quadrature that approximates the univariate integral%
\begin{equation*}
\int\limits_{a}^{b}G\left( x\right) dx
\end{equation*}%
is given by%
\begin{equation*}
Q_{C}\left[ G\left( x\right) \right] =\sum\limits_{i=1}^{N}c_{i}G\left(
x_{i}\right) =h\sum\limits_{i=1}^{N}w_{i}G\left( x_{i}\right) ,
\end{equation*}%
where the $x_{i}$ are nodes on $\left[ a,b\right] ,$ the coefficients $c_{i}$
are appropriate \textit{weights,} $h$ is a stepsize parameter representing
the separation of the nodes, and the \textit{reduced} weights are $%
w_{i}=c_{i}/h.$

The bivariate integral%
\begin{equation*}
\int\limits_{a}^{b}\int\limits_{l\left( x\right) }^{u\left( x\right)
}G\left( x,y\right) dydx
\end{equation*}%
is approximated by%
\begin{equation}
Q_{C}\left[ G\left( x,y\right) \right] =h\sum\limits_{i=1}^{N_{1}}w_{i}%
\left( k_{i}\sum\limits_{j=1}^{N_{2,i}}v_{j,i}G\left( x_{i},y_{j,i}\right)
\right) ,  \label{Qc[G(x,y)] = h wi ki vji}
\end{equation}%
where $v_{j,i}$ are appropriate reduced weights, $y_{j,i}$ are nodes along
the $y$-axis on $\left[ l\left( x_{i}\right) ,u\left( x_{i}\right) \right] ,$
and $k_{i}$ are stepsizes, with%
\begin{equation*}
k_{i}=\frac{u\left( x_{i}\right) -l\left( x_{i}\right) }{N_{2,i}}.
\end{equation*}%
Clearly, bivariate cubature is based on univariate quadrature. We can write%
\begin{equation}
Q_{C}\left[ G\left( x,y\right) \right] =\sum\limits_{i=1}^{N_{1}}\sum%
\limits_{j=1}^{N_{2,i}}C_{j,i}G\left( x_{i},y_{j,i}\right) ,
\label{Qc[G(x,y)] = double sum N1 N2,i}
\end{equation}%
where%
\begin{equation*}
C_{j,i}\equiv hw_{i}k_{i}v_{j,i}.
\end{equation*}

\subsection{Approximation error}

The approximation error in $Q_{C}\left[ G\left( x\right) \right] $ is
bounded by%
\begin{equation*}
A\left( r\right) \left( b-a\right) h^{r}\max_{\left[ a,b\right] }\left\vert
G^{\left( r\right) }\right\vert ,
\end{equation*}%
where $A\left( r\right) $ is a numerical constant particular to the type of
quadrature used (e.g. Trapezium, Simpson, Gauss-Legendre), and $r$ indicates
the so-called \textit{order} of the quadrature. Hence, for bivariate
cubature we have 
\begin{equation*}
A\left( r\right) \left( b-a\right) \underset{D}{\underbrace{\max \left(
u\left( x\right) -l\left( x\right) \right) }}\left( h^{r}\max \left\vert 
\frac{\partial ^{r}G}{\partial x^{r}}\right\vert +\left( \max
k_{i}^{r}\right) \max \left\vert \frac{\partial ^{r}G}{\partial y^{r}}%
\right\vert \right)
\end{equation*}%
as an upper bound on the approximation error. The integers $N_{1}$ and $%
N_{2,i}$ in (\ref{Qc[G(x,y)] = double sum N1 N2,i}) can be determined by
setting $h=\max k_{i}$ in the above bound, and demanding%
\begin{align}
& h^{r}A\left( r\right) \left( b-a\right) D\left( \max \left\vert \frac{%
\partial ^{r}G}{\partial x^{r}}\right\vert +\max \left\vert \frac{\partial
^{r}G}{\partial y^{r}}\right\vert \right) \leqslant \varepsilon  \notag \\
& \Rightarrow h=\left( \frac{\varepsilon }{A\left( r\right) \left(
b-a\right) D\left( \max \left\vert \frac{\partial ^{r}G}{\partial x^{r}}%
\right\vert +\max \left\vert \frac{\partial ^{r}G}{\partial y^{r}}%
\right\vert \right) }\right) ^{\frac{1}{r}},  \label{h = (...)^(1/r)}
\end{align}%
where the various maxima are found on the region of integration. Then%
\begin{align}
N_{1}& =\left\lceil \frac{b-a}{h}\right\rceil  \notag \\
N_{2,i}& =\left\lceil \frac{u\left( x_{i}\right) -l\left( x_{i}\right) }{k}%
\right\rceil .  \label{N2,i}
\end{align}%
Furthermore, the stepsizes $h$ and $k$ must be recalculated to be consistent
with integer values of $N_{1}$ and $N_{2,i},$ as in%
\begin{align}
h^{\ast }& =\frac{b-a}{N_{1}}  \notag \\
k_{i}^{\ast }& =\frac{u\left( x_{i}\right) -l\left( x_{i}\right) }{N_{2,i}},
\label{k*i}
\end{align}%
and it is these stepsizes that are used in (\ref{Qc[G(x,y)] = h wi ki vji}).
Once the stepsizes have been determined, the nodes $x_{i}$ and $y_{j,i}$ can
be found.

This process of computing stepsizes consistent with a tolerance $\varepsilon 
$ constitute \textit{absolute} error control in bivariate composite
interpolatory cubature, and is used in the previously described algorithm to
find $Q_{C}\left[ g\left( w,z\right) \right] $ such that%
\begin{equation*}
\left\vert I\left[ g\left( w,z\right) \right] -Q_{C}\left[ g\left(
w,z\right) \right] \right\vert \leqslant \varepsilon .
\end{equation*}

It should be noted that our use of $\max \left\vert \frac{\partial ^{r}G}{%
\partial x^{r}}\right\vert +\max \left\vert \frac{\partial ^{r}G}{\partial
y^{r}}\right\vert $ is conservative, and could result in smaller stepsizes
than is necessary, for the given tolerance. However, in these types of
numerical calculations it is always better to err on the side of caution.
Nevertheless, we should be aware that such a conservative approach could
result in $\left\vert I\left[ g\left( w,z\right) \right] -Q_{C}\left[
g\left( w,z\right) \right] \right\vert \ll \varepsilon ,$ so that $\frac{%
\varepsilon }{\left\vert Q_{C}\left[ g\left( w,z\right) \right] \right\vert }
$ overestimates the relative error. Analytically speaking, the approximation
error is proportional to%
\begin{equation*}
\left. \frac{\partial ^{r}G}{\partial x^{r}}\right\vert _{\left( \xi ,\zeta
\right) }+\left. \frac{\partial ^{r}G}{\partial y^{r}}\right\vert _{\left(
\varphi ,\phi \right) },
\end{equation*}%
where $\left( \xi ,\zeta \right) $ and $\left( \varphi ,\phi \right) $ are
points somewhere in the region of integration - but since these points are
not known, and we cannot be sure of the sign of the derivatives, we use $%
\max \left\vert \frac{\partial ^{r}G}{\partial x^{r}}\right\vert +\max
\left\vert \frac{\partial ^{r}G}{\partial y^{r}}\right\vert $ in the error
term, instead.

\subsection{Choosing between absolute and relative error control}

From (\ref{rel err < eps}) we have%
\begin{equation*}
\left\vert I\left[ G\left( x,y\right) \right] -Q_{C}\left[ G\left(
x,y\right) \right] \right\vert \leqslant \varepsilon \left\vert I\left[
G\left( x,y\right) \right] \right\vert ,
\end{equation*}%
so that relative error control is equivalent to absolute error control with
an effective tolerance $\varepsilon \left\vert I\left[ G\left( x,y\right) %
\right] \right\vert .$ Replacing $\varepsilon $ in (\ref{h = (...)^(1/r)})
with $\varepsilon \left\vert I\left[ G\left( x,y\right) \right] \right\vert $
shows that, if $\left\vert I\left[ G\left( x,y\right) \right] \right\vert
>1, $ $h$ would be larger than if the tolerance was simply $\varepsilon ,$
and if $\left\vert I\left[ G\left( x,y\right) \right] \right\vert <1,$ $h$
would be smaller. Consequently, $N_{1}$ and $N_{2,i}$ would be smaller or
larger, respectively. Smaller values of $N_{1}$ and $N_{2,i}$ imply greater
computational efficiency and so, for the sake of efficiency, we choose
relative error control when $\left\vert I\left[ G\left( x,y\right) \right]
\right\vert >1,$ and absolute error control when $\left\vert I\left[ G\left(
x,y\right) \right] \right\vert <1.$ When $\left\vert I\left[ G\left(
x,y\right) \right] \right\vert =1,$ the two cases are identical. This is why
we can impose absolute error control on $\left\vert I\left[ g\left(
w,z\right) \right] -Q_{C}\left[ g\left( w,z\right) \right] \right\vert $ -
by our definition of $g$, $I\left[ g\left( w,z\right) \right] $ is
guaranteed to have a magnitude less than or equal to one.

\subsection{Roundoff error}

It is easily shown (see\ Appendix) that the roundoff error associated with (%
\ref{Qc[G(x,y)] = h wi ki vji}) is bounded by%
\begin{equation*}
4\left( b-a\right) D\mu ,
\end{equation*}%
where $\mu $ is a bound on the magnitude of the machine precision of the
finite-precision computing device being used, $\left\vert G\left( x,y\right)
\right\vert \leqslant 1$ on the region of integration, and the cubature used
is based on \textit{positive-weight} quadrature. In such quadrature, all
weights are positive; examples of such quadrature include the Trapezium
rule, Simpson's rule and all types of Gaussian quadrature. If the region of
integration has unit area, as does $I\left[ g\left( w,z\right) \right] $,
then the roundoff error simply has the bound $4\mu .$ The roundoff error
represents the minimum achievable accuracy in the cubature approximation,
and is incorporated into the error control procedure by replacing the
numerator of (\ref{h = (...)^(1/r)}) with%
\begin{equation*}
\varepsilon -4\left( b-a\right) D\mu .
\end{equation*}%
Clearly, it makes no sense to impose a tolerance smaller than the roundoff
error. A typical desktop PC has $\mu \sim 10^{-16}.$

\section{Numerical examples}

\subsection{Example I: relative error control}

We approximate%
\begin{equation*}
I\left[ G\left( x,y\right) \right] =\int\limits_{1}^{2}\int%
\limits_{x^{2}/5}^{x^{3}/5}e^{4xy}dydx=1.92660\times 10^{3}
\end{equation*}%
using Simpson's rule $\left( r=4,A\left( r\right) =\frac{16}{180}\right) $.
For ease of presentation we show all numerical values truncated to five
decimals or fewer, although all our calculations are performed in double
precision. The application of the algorithm to this example will be
described in detail. With the transformations (using $u_{1}=8/5,l_{1}=1)$ 
\begin{align*}
x& =w+1,y=\frac{7z}{5}+\frac{1}{5} \\
& \left( \Rightarrow m_{1}=1,m_{2}=\frac{7}{5}\right) ,
\end{align*}%
the integral becomes%
\begin{align*}
I\left[ G\left( w,z\right) \right] & =\int\limits_{0}^{1}\int\limits_{%
\widetilde{l}\left( w\right) }^{\widetilde{u}\left( w\right) }\frac{7}{5}%
e^{4\left( w+1\right) \left( \frac{7z}{5}+\frac{1}{5}\right) }dzdw \\
\widetilde{u}\left( w\right) & =\frac{7\left( w+1\right) ^{3}}{25}-\frac{7}{5%
} \\
\widetilde{l}\left( w\right) & =\frac{7\left( w+1\right) ^{2}}{25}-\frac{7}{5%
}.
\end{align*}%
We find%
\begin{align*}
M& =5.07104\times 10^{5} \\
\max \left\vert \frac{\partial ^{4}g}{\partial w^{4}}\right\vert &
=1.67772\times 10^{3} \\
\max \left\vert \frac{\partial ^{4}g}{\partial z^{4}}\right\vert &
=1.57351\times 10^{4} \\
D& =\max \left( \widetilde{u}\left( w\right) -\widetilde{l}\left( w\right)
\right) =\frac{18}{26}.
\end{align*}%
The stepsize $h$ is given by%
\begin{equation}
h=\left( \frac{\varepsilon -4\mu }{\left( \frac{16}{180}\right) \left(
1\right) \left( \frac{18}{26}\right) \left( \max \left\vert \frac{\partial
^{4}g}{\partial w^{4}}\right\vert +\max \left\vert \frac{\partial ^{4}g}{%
\partial z^{4}}\right\vert \right) }\right) ^{\frac{1}{4}}=5.52707\times
10^{-4}  \label{h = ....00356}
\end{equation}%
and so, with $\varepsilon =10^{-10},$%
\begin{equation*}
N_{1}=1810,h^{\ast }=5.52486\times 10^{-4}
\end{equation*}%
Here, $h^{\ast }$ is the length of each simpson subinterval (which contains
three nodes), and there are $1810$ such subintervals. Hence, there are $3621$
nodes $w_{i}$ on $\left[ 0,1\right] $ with separation $h^{\ast }/2$ (this is
the reason for the factor $16=2^{4}$ in (\ref{h = ....00356})).

The stepsizes $k_{i}^{\ast }$ along the $z$-axis are found from (\ref{N2,i})
and (\ref{k*i}) for each $i=1,2,\ldots ,563,$ and we find%
\begin{equation*}
\max k_{i}^{\ast }=5.52706\times 10^{-4}.
\end{equation*}%
This enables the nodes $z_{j,i}$ $\left( j=1,2,\ldots ,N_{2,i}\right) $ to
be found, for each $i$. As with $w_{i}$, the spacing between these nodes is $%
k_{i}^{\ast }/2.$ It must be noted that $N_{2,i}$ could be zero, in which
case $k_{i}^{\ast }$ will be NaN (\textit{not-a-number} in IEEE arithmetic).
In such cases, it is appropriate to simply set $k_{i}^{\ast }=0.$

Composite Simpson quadrature of $g\left( w,z\right) $ is performed along the 
$z$-axis, for each $i,$ yielding the $3621$ quantities $Q_{C}\left[ g\left(
w_{i},z\right) \right] ,$ which have the form 
\begin{equation*}
Q_{C}\left[ g\left( w_{i},z\right) \right] =\frac{k_{i}^{\ast }}{6}\left[ 
\begin{array}{c}
g\left( w_{i},z_{1,i}\right) +4g\left( w_{i},z_{2,i}\right) +2g\left(
w_{i},z_{2,i}\right) +4g\left( w_{i},z_{4,i}\right) + \\ 
\ldots +2g\left( w_{i},z_{N_{2,i}-2,i}\right) +4g\left(
w_{i},z_{N_{2,i}-1,i}\right) +g\left( w_{i},z_{N_{2,i},i}\right)%
\end{array}%
\right] .
\end{equation*}%
The integer coefficients in this expression are the weights appropriate to
composite Simspon quadrature.

Finally, Simpson quadrature is performed over these quantities along the $w$%
-axis, to give%
\begin{align*}
Q_{C}\left[ g\left( w,z\right) \right] & =\frac{h_{i}^{\ast }}{6}\left[ 
\begin{array}{c}
Q_{C}\left[ g\left( w_{1},z\right) \right] +4Q_{C}\left[ g\left(
w_{2},z\right) \right] +2Q_{C}\left[ g\left( w_{3},z\right) \right] + \\ 
4Q_{C}\left[ g\left( w_{4},z\right) \right] +\ldots +2Q_{C}\left[ g\left(
w_{N_{1}-2},z\right) \right] + \\ 
4Q_{C}\left[ g\left( w_{N_{1}-1},z\right) \right] +Q_{C}\left[ g\left(
w_{N_{1}},z\right) \right]%
\end{array}%
\right] \\
& =3.79922\times 10^{-3}.
\end{align*}

Hence, 
\begin{equation*}
I\left[ G\left( x,y\right) \right] \approx MQ_{C}\left[ g\left( w,z\right) %
\right] =1.92660\times 10^{3}.
\end{equation*}%
The estimate of the relative error is%
\begin{equation*}
\left\vert \frac{\varepsilon }{Q_{C}\left[ g\left( w,z\right) \right] }%
\right\vert =2.63211\times 10^{-8}
\end{equation*}%
while the actual relative error is $1.47276\times 10^{-11}.$ Clearly, the
actual error is less than the estimate. This is to be expected when using $%
\max \left\vert \frac{\partial ^{r}G}{\partial x^{r}}\right\vert +\max
\left\vert \frac{\partial ^{r}G}{\partial y^{r}}\right\vert $ in the
computation of $h.$ Obviously, our value for $h$ is conservative (smaller
than actually necessary) and so the actual error is smaller than the
estimate. Nevertheless, as we have stated earlier, it is better to be more
accurate than necessary and, since the estimate reflects an upper bound, we
can be sure that the error is no more than $2.63211\times 10^{-8}.$ If this
level of accuracy is acceptable, then the result stands. However, if we
desire a relative error of no more than $10^{-10},$ say, we simply repeat
the algorithm with%
\begin{equation*}
\frac{\varepsilon }{264}
\end{equation*}%
as the new tolerance. This gives%
\begin{equation*}
\left\vert \frac{\varepsilon }{Q_{C}\left[ g\left( w,z\right) \right] }%
\right\vert =9.97014\times 10^{-11}<10^{-10},
\end{equation*}%
while the actual relative error is $5.35801\times 10^{-14}.$

\subsection{Example II: absolute error control}

In this second example, the integral 
\begin{equation*}
I\left[ G\left( x,y\right) \right] =\int\limits_{1}^{4}\int%
\limits_{x}^{2x^{2}}\frac{\sin \left( xy\right) }{5}dydx=-0.00734
\end{equation*}%
will again be approximated using Simpson quadrature but, since it has
magnitude less than one, we will see that absolute error control is more
efficient than relative error control. There is no need for a detailed
exposition, as in the previous example, and we simply state our results.

The transformed integral is%
\begin{align*}
I\left[ G\left( w,z\right) \right] & =\int\limits_{0}^{1}\int\limits_{%
\widetilde{l}\left( w\right) }^{\widetilde{u}\left( w\right) }\frac{93}{5}%
\sin \left( \left( 3w+1\right) \left( 31z+1\right) \right) dzdw \\
\widetilde{u}\left( w\right) & =\frac{2\left( 3w+1\right) ^{2}-1}{31} \\
\widetilde{l}\left( w\right) & =\frac{3w}{31}.
\end{align*}%
Using $M=\frac{93}{5}$ and $\varepsilon =10^{-4}$ gives%
\begin{align*}
M\varepsilon & =0.00186 \\
\left\vert \frac{\varepsilon }{Q_{C}\left[ g\left( w,z\right) \right] }%
\right\vert & =0.25340.
\end{align*}%
The upper bound on the relative error is fairly large. The absolute error is
estimated by $M\varepsilon $ ; it is clear that, since $M$ is known, $%
\varepsilon $ can be chosen so that $M\varepsilon $ equals some desired
value. For example, if we seek an absolute error of no more than $10^{-5},$
we choose $\varepsilon =5.37\times 10^{-7},$ which gives%
\begin{align*}
M\varepsilon & =9.9882\times 10^{-6}<10^{-5} \\
\left\vert \frac{\varepsilon }{Q_{C}\left[ g\left( w,z\right) \right] }%
\right\vert & =0.00136,
\end{align*}%
with $N_{1}=1761$ (hence, $3523$ nodes on the $w$-axis). Note that achieving
this tolerance does not require a repetition of the algorithm, since $M$ is
known a priori.

On the other hand, to improve the estimate of the relative error to $10^{-5}$
requires $\varepsilon =10^{-5}/136=7.353\times 10^{-8},$ which results in $%
N_{1}=2895,$ and hence, more nodes than are needed to achieve the same
tolerance in the absolute error. This is consistent with our earlier
discussion regarding the efficiency-based criterion for choosing between
absolute and relative error control.

\section{Conclusion}

We have reported on an algorithm for controlling the relative error in the
numerical approximation of a bivariate integral. The numerical method used
is positive-coefficient composite interpolatory quadrature. The algorithm
involves transforming and scaling the integral to one that has magnitude
bounded by unity, and then applying an absolute error control procedure to
such integral. The relevant scaling factor is then used to find the
approximate value of the original integral and an estimate of the relative
error (if the integral has magnitude greater than unity) or absolute error
(if the integral has magnitude less than or equal to unity). The calculation
can be repeated with an appropriate refinement if the estimated error is
considered too large. The algorithm proceeds in a systematic and controlled
manner, and there is no need for any prior knowledge of the magnitude of the
integral. Two examples with Simpson's rule clearly demonstarte the character
of the algorithm. This work extends other work of ours \cite{jscp}, in which
we considered the control of relative error in the quadrature of a
univariate integral. In that work, we designated the algorithm CIRQUE, and
so we take the liberty here of designating the current algorithm CIRQUE2D.

\begin{flushleft}
\bigskip \textbf{{\Large {Disclosure statement}}}\ \ 

The author reports there are no competing interests to declare.

\begin{equation*}
\end{equation*}%
\textbf{{\Large {Appendix:}}}\ \ \textbf{Roundoff bound}
\end{flushleft}

Using (\ref{Qc[G(x,y)] = h wi ki vji}) an (\ref{Qc[G(x,y)] = double sum N1
N2,i}), we write%
\begin{align*}
Q_{C}\left[ G\left( x,y\right) \right] &
=\sum\limits_{i=1}^{N_{1}}c_{i}\left( 1+\mu _{c,i}\right) \\
& \times \sum\limits_{j=1}^{N_{2,i}}k_{i}\left( 1+\mu _{k,i}\right)
v_{j,i}\left( 1+\mu _{v,j,i}\right) G\left( x_{i},y_{j,i}\right) \left(
1+\mu _{G,j,i}\right) \\
& =\sum\limits_{i=1}^{N_{1}}\sum\limits_{j=1}^{N_{2,i}}C_{j,i}G\left(
x_{i},y_{j,i}\right) +C_{j,i}G\left( x_{i},y_{j,i}\right) \left( \mu
_{w,i}+\mu _{v,j,i}+\mu _{G,j,i}\right) ,
\end{align*}%
where we have indicated the roundoff error in $c_{i},k_{i},v_{j,i}$ and $%
G\left( x_{i},y_{j,i}\right) $ explicitly, and we have ignored higher-order
terms in the second line. The roundoff error $\Upsilon $ in the cubature is%
\begin{align*}
\Upsilon & \equiv
\sum\limits_{i=1}^{N_{1}}\sum\limits_{j=1}^{N_{2,i}}C_{j,i}G\left(
x_{i},y_{j,i}\right) \left( \mu _{c,i}+\mu _{w,i}+\mu _{v,j,i}+\mu
_{G,j,i}\right) \\
& \leqslant \sum\limits_{i=1}^{N_{1}}\sum\limits_{j=1}^{N_{2,i}}4C_{j,i}\mu ,
\end{align*}%
where $\mu $ is a bound on $\left\vert \mu _{c,i}\right\vert ,\left\vert \mu
_{w,i}\right\vert ,\left\vert \mu _{v,j,i}\right\vert $ and $\left\vert \mu
_{G,j,i}\right\vert ,$ and we have assumed $\left\vert G\left(
x_{i},y_{j,i}\right) \right\vert \leqslant 1.$ Now, since $%
C_{j,i}=hw_{i}k_{i}v_{j,i}=c_{i}k_{i}v_{j,i},$%
\begin{equation*}
\Upsilon \leqslant 4\mu \sum\limits_{i=1}^{N_{1}}c_{i}\left(
\sum\limits_{j=1}^{N_{2,i}}k_{i}v_{j,i}\right) .
\end{equation*}%
But, in positive-weight univariate composite interpolatory quadrature, the
sum of the weights is simply the length of the interval of integration, and
so%
\begin{align*}
\Upsilon & \leqslant 4\mu \left( b-a\right) \left( \max \left( u\left(
x_{i}\right) -l\left( x_{i}\right) \right) \right) \\
& =4\mu \left( b-a\right) D \\
& \leqslant 4\mu
\end{align*}%
if $\left( b-a\right) D\leqslant 1.$

\end{document}